%Paper: alg-geom/9502004
%From: masanori@mathp7.jussieu.fr (Masanori Kobayashi)
%Date: Mon, 6 Feb 95 14:41:10 +0100

\input amstex		%AmSTeX 2.1
\documentstyle{amsppt}
\magnification=\magstep1
 
%following symbols can be used in the math mode.
\define\bb{\Bbb}	%if you don't have Blackboard bold, change here
 \define\C{{\bb C}}   \define\R{{\bb R}}
\define\Q{{\bb Q}}   \define\Z{{\bb Z}}   \define\N{{\bb N}}
\redefine\P{{\bb P}}    \define\T{{\bb T}}
 \define\CO{{\Cal O}} 
\define\Proj{\operatorname{Proj}}

\define\Hom{\text{\rm Hom}}       \define\Aut{\text{\rm Aut}}
\define\rk{\operatorname{rk}}
\define\Int{\operatorname{Int}}

\define\chap#1{\medskip\noindent{\bf #1}\nopagebreak\par}

\topmatter
 \title Duality of weights, mirror symmetry and Arnold's strange
duality\endtitle
 \author Masanori KOBAYASHI\endauthor
 \affil Department of Mathematics, Tokyo Institute of Technology and
  U.F.R de Math\'ematique, Universit\'e Paris VII\endaffil
 \address
  Department of Mathematics, Tokyo Institute of Technology, 2-12-1 Oh-okayama,
Meguro, Tokyo 152 Japan;
  U.F.R de Math\'ematique, Universit\'e Paris VII, 45-55, 5\`eme \'etage, 2
Place Jussieu, 75251 Paris cedex 05, France
 \endaddress
 \email masanori\@math.titech.ac.jp, masanori\@mathp7.jussieu.fr\endemail
 \thanks
  Partially supported by F\^ujukai Foundation
  and Japan Association for Mathematical Sciences.
 \endthanks
\endtopmatter
\document
%\baselineskip=2\baselineskip	%please do not comment out at last.

\abstract
A notion of duality of weight systems which corresponds to Batyrev's toric
mirror symmetry is given. Explicit duality on the (1,1)-cohomology of K3
surfaces which are minimal models of toric hypersurfaces is constructed using
monomial divisor mirror map of Aspinwall-Greene-Morrison. It is shown that
Arnold's strange duality for exceptional unimodal singularities reduces to
this duality.
\endabstract

\chap{Introduction.}
The hypersurfaces in weighted projective spaces often appear
as important examples in the context of mirror symmetry.
In this paper, we describe the relation between polar duality
and duality of weight systems.

The duality of weights partly suggests why \cite{CLS}
produced a mirror symmetric phenomena
using only a resolution of weighted hypersurface in weighted $\P^4$.
In fact, recently it is shown that those examples in weighted 4-spaces
correspond to some reflexive polytopes \cite{CdOK}.

As an application, we will show that Arnold's strange duality
for fourteen unimodal singularities reduces to polar duality.
\medskip

The relation between the weight systems and
the hypersurface singularities with a $\C^*$-action is as follows;
take a germ of analytic function $f:\C^n \to \C$ with $f(0)=0$, which
determines a germ of hypersurface singularity $(\{f=0\},0) \subset (\C^n,0)$.
The function $f$ can be expressed as a weighted homogeneous polynomial
of $n$ variables $(x_1,\dots,x_n)$ after some suitable analytic change of
coordinates near the origin if and only if
there is a holomorphic tangent vector field $D$ such that $Df=f$\cite{S1}.
In this case, we can assign to $f$
a weight system $W=(wt(x_1),\dots,wt(x_n);wt(f))$.

We associate some hypersurface singularity with $\C^*$-action
a compact complex surface with trivial canonical sheaf
in the following manner.

Let $X_0 \subset \C^n$ be a hypersurface singularity
with $\C^*$-action whose weights are all positive.
Such $X_0$ is known to be an algebraic variety defined
by a weighted homogeneous polynomial\cite{OW}.

Let $X_t$ be its Milnor fiber (if it exists;
e.g. when $X_0$ is an isolated singularity).
We will not restrict ourselves to isolated singularities.
Let $\bar X$ be the natural compactification of $X_t$ in $\P(1,a_1,\dots,a_n)$
by sending the points $(x_1,\dots,x_n)$ of $\C^n$ to $(1:x_1:\dots:x_n)$.

We treat only the case that total degree $h$ is greater than the sum of
the weights $\sum_{i=1}^n a_i$ here.
We denote the difference $h-\sum_{i=1}^n a_i$ by $a_0$.

Let $\hat X$ be the image of $\bar X$ by the $\Z/a_0\Z$-quotient:
$\P(1,a_1,\dots,a_n)\to\P(a_0,a_1,\dots,a_n)$.
By the adjunction formula, the canonical sheaf of $\hat X$ is trivial.
We denote a crepant resolution of $\hat X$ by $X$ if it exists,
in which case $\hat X$ has only canonical Gorenstein singularities.

In the case $X_0$ is an isolated singularity,
it is known that the Milnor lattice or exceptional lattice of resolution
and the exceptional lattice at infinity enjoy some dual property\cite{LW}.

We assume $X$ is a K3 surface corresponding to one of Arnold's
singularities\cite{P}\cite{DN}.

The transcendental lattice $T_X$ decomposes as $L_G \oplus U$,
where the dual graph of $L_G$ is a tree of $(-2)$-elements with three
branches and $U$ is the even unimodular hyperbolic lattice of rank $2$.
We will write $L_D := S_X$ in this case.
Then the whole cohomology ring decomposes as
$H^*(X,\C)=H^{0,0}\oplus (L_G \oplus U \oplus L_D)_\C \oplus H^{2,2}$,
where the subscript ${}_\C$ means the tensor product by $\C$.
The mirror map should interchange $(L_G)_\C$ and $(L_D)_\C$,
$H^{0,0}\oplus H^{2,2}$ and $U_\C$, respectively,
since the vector space $(L_G)_\C$ corresponds to the tangent space of
the deformation space of complex structures with a fixed Picard number
and $(L_D)_\C$ corresponds to a deformation of
complexified algebraic K\"ahler structures.
This is an explanation due to \cite{AM}
why mirror symmetry for algebraic K3 surface
is said to correspond to Arnold's duality.

Roan\cite{Ro} applied orbifold construction for four of the dual pairs
and constructed non-linear coordinate changes for them.
Borcea\cite{Bo} used Batyrev's polar dual construction
for the reflexive pyramids and calculated some of dual weights.

We shall give a proof that polar dual of
polytopes corresponding Arnold's singularities
gives corresponding dual singularity.

More generally, take a K3 surface $X$ which is the minimal resolution
of toric hypersurface with only rational double points.
Let $L_D$ be the restriction of the toric divisors of the ambient space
and $L_0$ be the orthogonal complement $L_D^{\bot}$ in $S_X$.
then $(S_X)_\Q$ decomposes as a direct sum $(L_D)_\Q\oplus (L_0)_\Q$.

We construct the monomial divisor mirror map \cite{AGM} $\mu$ in this context
which interchanges $(L_D)_\C$ and $(L_G)_\C$, and fixes $(L_0)_\C$,
and show that $\mu$ corresponds to polar duality.

The reason why the full algebraic lattice and $L_G$ does not
interchange is that the induced complexified K\"ahler class is trivial
on $L_0$-part while embedded deformation has the full dimension.
This $L_0$-part is a generalization of $H^{1,1}_n$ in \cite{Ro},
where Roan showed a similar mirror property for a quotient construction.
Since the mirror of generic marked K3 surface seems to be again a generic one,
mirror phenomena for special marked K3 surfaces are of particular interest,
where `special' means its (complexified) K\"ahler class is also special.

For general K3 surfaces
$L_0$ is defined as the orthogonal complement of the fixed lattice by
the automorphism group which fixes $H^{2,0}$ and the K\"ahler class.
We show this group is finite
using Nikulin's results on the automorphism of a K3 surface and Torelli's
theorem.

We refer to \cite{AM}\cite{K} and their reference for an introduction to some
background and \cite{Mr}\cite{Y} for general information about mirror symmetry.
Very recently, the author receives a preprint by Dolgachev \cite{D2}
which includes relevant general results and examples using lattice theory.

In Chapter 1, we shall formulate a dual correspondence of reflexive pairs
in $\Q$-Fano toric projective varieties
after reviewing some results of toric geometry.
In Chapter 2, we prepare the language of weight systems according to K. Saito
and introduce our notion of duality of weight systems.
Our duality partly coincides with Saito's duality.
In Chapter 3, we state the relation between duality of weight systems and
Batyrev's toric mirror symmetry using polar polytopes.
We found that our construction is partly a generalization of those of
Batyrev\cite{B1}, Borcea\cite{Bo} and Berglund-Katz\cite{BK}
for reflexive simplices and pyramids.
In Chapter 4, we review some general results of K3 surfaces,
and we argue in the case of K3 surfaces and construct a monomial
divisor mirror map of Aspinwall-Greene-Morrison and prove
Arnold's strange duality.
We add some examples which are not Arnold's singularities.

\chap{Acknowledgement}

The author expresses his hearty gratitude to Y. Kawamata
for his continuous warm encouragement,
F. Campana, D.-T. L\^e, H. Popp, M. Reid, B. Teissier,
C. Weber and P. M. H. Wilson
for their excellent hospitality, stimulus discussion and valuable suggestions.

He also thanks
to M. Furushima, F. Michele, I. Nakamura and J. Stevens
for useful knowledge of singularity theory,
to K. Saito for a very interesting introduction of
his theory of weight systems,
to V. V. Batyrev, D. R. Morrison, Y. Ruan and A. Tyurin for useful information
and to H. Ooguri, Y. Shimizu for warm encouragement.

\chap{Notation.}
$\N$ : the monoid of nonnegative integers.\par
$\Z_+$ : the semigroup of positive integers.\par
$\Z$ : the ring of integers.\par
$\Q,\R,\C$ : the field of rational, real, complex numbers.\par
$\T$ : $n$-dimensional complex torus $(\C^*)^n$. \par

\chap{1. Polar duality.}
\chap{1.1 Reflexive polytopes.}
We fix here some notations about geometry of convex body and
review some results.
We refer \cite{O} for generality of
toric geometry and geometry of convex bodies.

Let $M$ be a free $\Z$-module of rank $n$.
We naturally imbed $M$ in the $\Q$-vector space $M_\Q:=M\otimes_\Z\Q$.
Similarly for the dual module
$N:=M^\vee=\Hom(M,\Z)$.
We denote the natural pairing by $\left<\,,\,\right>:M\times N \to \Z$,
which we also use for its $\Q$-extensions.

\definition{Definition 1.1.1}
A subset $\Delta$ in $M_\Q$ is called a {\it convex polytope\/}
if it is a convex hull of a finite subset of $M_\Q$.
$\Delta$ is {\it integral\/} if the finite subset can be taken from $M$.
A {\it face\/} of $\Delta$ is a nonempty intersection with a hyperplane whose
closed half space contains whole $\Delta$.
We denote the {\it $k$-skeleton\/} of $\Delta$ by $\Delta^{[k]}$,
which is the union of all faces of $\Delta$
whose dimension is no more than $k$.
We call a codimension-one face a {\it facet\/}.
We denote $|M\cap\Delta|$ and $|M \cap\Int\Delta|$ by
$l(\Delta)$ and $l^*(\Delta)$, respectively.
\enddefinition

\remark{Remark 1.1.2}
We use those words above similarly for subsets in $N_\Q$.
We shall sometimes omit the word ``convex" in this paper
since we always assume so for polytopes.
A face of a polytope $\Delta$ is also a polytope.
\endremark

\definition{Definition 1.1.3}
Let $K$, $K_1$ and $K_2$ be subsets of $M_\Q$ and $c\in\Q$.
$cK := \{cx \in M_\Q | x \in K\}$,
$K_1 + K_2 := \{x_1 + x_2 \in M_\Q | x_i \in K_i \text{ for }i = 1,2\}$.
\enddefinition

\definition{Definition 1.1.4}
$K$ be a subset in $M_\Q$.
The {\it polar dual\/} of $K$ is the following subset of $N_\Q$:
$K^*$=$\{y\in N_\Q | \left<x,y\right>\ge-1 \text{ for all } x\in K \}$.
\enddefinition

\proclaim{Lemma 1.1.5}
\roster
\item $K^*$ is a convex set containing $0$.
If $K_1 \subset K_2$ then $K_1^* \supset K_2^*$.
\item If $\Delta$ is an $n$-dimensional convex polytope
with $0\in \Int\Delta$, then all the same for $\Delta^*$.
Moreover, $\Delta^{**} = \Delta$ by natural identification $M^{\vee\vee} = M$.
\endroster
\endproclaim

\definition{Definition 1.1.6}\cite{Da}
For an $n$-dimensional polytope $\Delta$ in $M_\Q$,
define an $n$-dimensional projective variety
$$\P_{\Delta,M} = \P_\Delta := \Proj \bigoplus_{k=0}^\infty \C\left<k\Delta\cap
M\right>T^k,$$
where $\C\left<k\Delta\cap M\right>$ denotes the $\C$-vector space generated
the elements of $M$ in $k\Delta$, and the multiplication comes from the
addition of $M$.
$(\P_\Delta,\CO(1))$ is the polarized variety
associated to the polytope $\Delta$.
\enddefinition

\definition{Definition 1.1.7}\cite{B1}
Let $\Delta$ be an $n$-dimensional integral polytope with $0 \in \Int \Delta$.
$\Delta$ is said to be {\it reflexive\/}
if one of the following equivalent conditions is satisfied:
\roster
\item $\Delta^*$ is also an integral polytope;
\item there exists a finite subset $\{y_1,\dots,y_k\}$ in $N$ such that
$\Delta = \{x\in M_\Q\,|\,\left<x,y_i\right>\ge-1\text{ for }1\le i \le k\}$;
\item for each facet $\delta$ of $\Delta$,
there exists an integral element $y \in N$ such that
$\delta \subset \{x \in M_\Q\,|\,\left<x,y\right>=-1\}$;
\item for each facet $\delta$, there are no points of $M$ between
the origin and the hyperplane containing $\delta$.
\endroster
\enddefinition

If $\Delta$ is reflexive, $\Delta^*$ is reflexive.
We cite here some results in \cite{B1}.

\proclaim{Theorem 1.1.8}\cite{B1 Theorem 4.1.9}
Let $\Delta$ be an $n$-dimensional integral polyhedron in $M_\Q$,
$\P_\Delta$ the corresponding $n$-dimensional projective toric variety,
and ${\Cal F}(\Delta)$ the family of projective $\Delta$-regular hypersurfaces
${\overline Z}_f$ in $\P_\Delta$.
Then the following conditions are equivalent:
\roster
\item the family ${\Cal F}(\Delta)$ of $\Delta$-regular hypersurfaces
in $\P_\Delta$ consists of Calabi-Yau varieties with canonical singularities;
\item the ample invertible sheaf $\CO_\Delta(1)$ on the toric variety
$\P_\Delta$ is anticanonical(i.e. $\P_\Delta$ is a toric Fano variety with
Gorenstein singularities);
\item $\Delta$ contains only one integral point $m_0$ in its interior,
and $(\Delta - m_0, M)$ is a reflexive pair.
\endroster
\endproclaim

\proclaim{Theorem 1.1.9}\cite{B1 Theorem 4.4.3}
For any reflexive polyhedron $\Delta$ of dimension $n \ge 4$,
the Hodge number $h^{n-2,1}({\widehat Z}_f)$ of a MPCP-desingularization of
a $\Delta$-regular Calabi-Yau hypersurface ${\overline Z}_f \subset \P_\Delta$
equals the Picard number $h^{1,1}({\widehat Z}_g)$ of a MPCP-desingularization
of a $\Delta^*$-regular projective Calabi-Yau hypersurface
${\overline Z}_g \subset \P_{\Delta^*}$ corresponding to the dual reflexive
polyhedron $\Delta^*$.
\endproclaim
In the above two theorems, ``polyhedron" means ``polytope" in our paper.
We reserve the word ``polyhedron" for a three dimensional polytope.
MPCP-desingularization is a maximal projective crepant partial
toric embedded resolution.

\chap{1.2. Dual correspondence for toric projective varieties.}

In this section, we shall treat a family of anticanonical
members in a $\Q$-Fano toric variety $\P_\Delta$
where $\Delta$ may not be integral.
We will think not only the quotient family corresponding to the same $\Delta$
but also specialization corresponding to a sub-Newton polytope.

\definition{Definition 1.2.1}
Let ${\Cal P}_n$ be the set of $n$-dimensional pairs $(\Delta, M)$
where $M$ is a free $\Z$-module of rank $n$ and
$\Delta$ is a $n$-dimensional polytope in $M_\Q$
which contains $0$ in its interior,
such that $\P_{\Delta,M}$ is $\Q$-Fano and $\Delta$ represents
anticanonical divisor.
\enddefinition

We usually assume that $n \ge 2$.
Let $(\overline\Delta, \overline M)$ and $(\overline\nabla, \overline N)$ be
two elements of ${\Cal P}_n$.
We denote the dual groups $(\overline M)^\vee$ by $\underline N$ and
$(\overline N)^\vee$ by $\underline M$.
We also denote the polar duals ${\overline\Delta}^*$ and ${\overline\nabla}^*$
by $\underline\nabla$ and $\underline\Delta$, respectively.

\definition{Definition 1.2.2}
A $\Q$-linear isomorphism
$\sigma\: {\underline M}_\Q \to {\overline M}_\Q$ such that
$\sigma(\underline\Delta) \subset \overline\Delta$,
is called a {\it quasi-correspondence\/} between
$(\overline\Delta, \overline M)$ and $(\overline\nabla, \overline N)$.
If $\sigma(\underline\Delta) \cap \overline M$ generates $\overline M$ and
$\underline\Delta$ has a nonempty intersection with
each facet of $\overline\Delta$,
we say $\sigma$ is a {\it correspondence\/}.

If $\sigma(\underline M) = \overline M$,
we say it {\it dual correspondence\/}.
\enddefinition

We shall identify ${\underline M}_\Q$ and ${\overline M}_\Q$,
${\overline N}_\Q$ and ${\underline N}_\Q$
by $\sigma$ and ${}^t\sigma$, respectively.
We shall sometimes omit $\sigma$.

When $(\overline\Delta,\overline M)$ and
$(\overline\nabla,\overline N)$ have a dual correspondence $\sigma$,
we will use $M = \overline M = \sigma(\underline M)$ and
$N = \overline N = {}^t\sigma(\underline N)$.

A basic observation is the following

\proclaim{Proposition 1.2.3}
Assume there exists a dual correspondence between two elements
$(\overline\Delta, M)$ and $(\overline\nabla, N)$ of ${\Cal P}_n$.
Let $(\Delta,M)$ be a reflexive pair such that
$\underline\Delta \subset \Delta \subset \overline\Delta$.

Then there is a natural birational map
$\Phi \: \P_{\overline\Delta} \dashrightarrow \P_\Delta$
coming from the inclusion between the polytopes.
A general anticanonical member of $\P_\Delta$ is birational
to the corresponding member of $\P_{\overline\Delta}$.

Similarly, from $\underline\nabla \subset \Delta^* \subset \overline\nabla$,
a general anticanonical member of $\P_{\Delta^*}$ is
birational to that of $\P_{\overline\nabla}$.
\endproclaim

\demo{Proof}
We denote the sublinear system of the anticanonical linear system of
$\P_{\overline\Delta}$ corresponding to the polytope $\Delta$
by $\Lambda(\Delta)$.

By definition, any nonempty nonzero subsystem $\Lambda$ is
free from base points on the $n$-dimensional orbit of $\T$.

Since the set $\Delta\cap M$ generates the whole $M$,
$\Lambda(\Delta)$ defines a birational map
$\overline\varphi\:\P_{\overline\Delta} \dashrightarrow \P^{l(\Delta)-1}$.
The image $V$ is the projective variety
$\Proj\bigoplus_{k=0}^\infty \left(\sum_{i=1}^k(\Delta\cap M)\right)$,
which is associated to the graded algebra
generated by $\Delta\cap M$.

Each facet of $\overline\Delta$ corresponds to an $(n-1)$-dimensional orbit
of $\T$ on $\P_{\overline\Delta}$,
which has nonempty intersection with $\Delta$,
since $\Delta$ contains $\underline\Delta$.
Thus $\Lambda(\Delta)$ is free from base components
and the birational transform of a hyperplane section of $V$
in $\P_{\overline\Delta}$ is a member of $\Lambda(\Delta)$.
In particular, a general member of $\Lambda(\Delta)$ is
an irreducible divisor if $n>1$.

On the other hand, the anticanonical complete linear system on $\P_\Delta$ is
free from base points, since $\Delta$ is an integral polytope,
and defines a morphism $\varphi$ to the same $V$,
which is birational by the same reason as above.
Take $\Phi := \varphi^{-1}\circ\overline\varphi$.
\qed
\enddemo

\remark{Remark 1.2.4}
Thus, in the situation above, if a Calabi-Yau variety $\bar Z_\Delta$ is
a specialization of $\bar Z_{\Delta'}$,
then $\bar Z_{\Delta^*}$ is a generalization of $\bar Z_{{\Delta'}^*}$.
When you specialize $\bar Z$, the singularity goes bad.
As a result, we have more K\"ahler moduli and less complex structure moduli
for the desingularization $\tilde Z_\Delta$,
and contrary for $\tilde Z_{\Delta^*}$.
We hope this relation holds also outside toric varieties.
\endremark

\remark{Remark 1.2.5}
$\P_\Delta$ is isomorphic to $V$ if $n \le 3$.
Note that we have unique minimal models in this case.
\endremark

\proclaim{Proposition 1.2.6}
Let $\Delta$ be a reflexive polytope of dimension at most three,
and $(\P_\Delta,\CO(1))$ be the polarized variety associated to it.
Then $\CO(1)$ is simply generated and
defines projective normal embedding in $\P^{l(\Delta)-1}$.
\endproclaim

\demo{Proof}
The proof follows from the following
\enddemo
\proclaim{Proposition 1.2.7}
Let $\Delta$ be a reflexive polytopes in $M_\Q$
whose dimension is at most three,
and $k$ be an positive integer.
Then for any integral point $P \in k\Delta$,
there exist integral points $P_i \in \Delta$ ($i=1,\dots,k$) with
$P = P_1 + \dots + P_k$.
\endproclaim
\demo{Proof}
We assume $n=3$, since other cases are similar and easier.
For each $2$-dimensional face $\delta$ of $\Delta$,
we have a subdivision by integral simplices (triangles) such that
\roster
\item each simplex is elementary, that is, it does not contain points of $M$
other than the vertices,
\item each simplex has common points with other simplices
only on its edges and
\item the union of all simplices covers $\delta$.
\endroster

We note that an elementary simplex of dimension two is regular.
That means, if we denote the vertices of the simplex by $Q_0$, $Q_1$ and $Q_2$
and the plane containing $\delta$ by $\alpha$,
$Q_1-Q_0$ and $Q_2-Q_0$ spans $\alpha \cap M$.
Let $C$ be the closed cone over the simplex with apex at $0$.
Then any integral point in $C$ can be written as an integral combination of
$Q_i$'s, since $\Delta$ is reflexive.
\enddemo

\remark{Remark 1.2.8}
In higher dimensions, there exist elementary non regular simplices.

In general for an ample divisor $A$ on an algebraic K3 surface,
it is known that $3A$ is very ample.
\endremark

\proclaim{Corollary 1.2.9}
For a general toric hypersurface $\bar Z$ of dimension $1$ or $2$,
the induced morphism
$\bar Z\to \P^{l(\Delta)-2}$.
is a projectively normal embedding.
\endproclaim

\demo{Proof}
$H^1(\P_\Delta, \CO(k))$ vanishes for $k\ge0$.
\qed
\enddemo

In Chapter 3,
we shall apply these notions to the case of weighted projective spaces
with one particular choice of a homogeneous coordinate.

\chap{2. Duality for weight systems.}
\chap{2.1. Weighted projective space.}

We fix here some notations.
For general properties about weighted projective spaces,
see for example, \cite{Mi}\cite{D1}\cite{F}.

\definition{Definition 2.1.1}
Let $a_0, \dots, a_n$ be fixed positive integers.
We denote by $\P(a_0,\dots,a_n)$ the projective variety
$\Proj\C[x_0,\dots,x_n]$ where the degree of $x_i$ is $a_i$.
We call it a {\it weighted projective space\/} of weight $(a_0,\dots,a_n)$.
\enddefinition

For a positive integer $k$, there is a natural isomorphism
$\P(a_0,\dots,a_n)\cong\P(ka_0,\dots,ka_n)$.

\definition{Definition 2.1.2}
A weight $(a_0, \dots, a_n)$ is {\it reduced\/}
if $\gcd(a_0,\dots,a_n)=1$.
\enddefinition

We assume the weight is reduced.
Let $k$ be $\gcd(a_1,\dots,a_n)$.
Then $\P(a_0,a_1,\dots,a_n)\cong\P(a_0,a_1/k,\dots,a_n/k)$,
which leads to the following
\definition{Definition 2.1.3}
A reduced weight $(a_0,\dots,a_n)$ is {\it well-formed\/}
if for all $i, \gcd(a_j)_{j\ne i}=1$.
\enddefinition

We shall usually treat only well-formed weights.

\chap{2.2. Weight systems.}

We refer \cite{S2} for general statements of weight systems,
especially for $n=3$.

\definition{Definition 2.2.1}
A ($n+1$)-uple of positive integers $W=(a_1,\dots,a_n;h)$ is called
a {\it system of weights\/} or simply a {\it weight system\/}.
We always assume that $h \in \sum_{i=1}^n a_i$.
We call the integers $a_i$ as {\it weights\/} of $W$ and
the last weight $h$ the {\it degree\/} of $W$.
\enddefinition

\definition{Definition 2.2.2}
$W$ is said to be {\it reduced\/} if $\gcd(a_1,\dots,a_n,h)=1$ holds.
%If $n=3$,
\enddefinition
\remark{Remark 2.2.3}
$W$ is reduced then $\gcd(a_1,\dots,a_n) = 1$, since $h \in \sum_{i=1}^n a_i$.
\endremark

\definition{Definition 2.2.4}
$W=(a_1,\dots,a_n;h)$ and $W'=(a'_1,\dots,a'_n;h')$ are {\it equivalent\/}
if for some rational number $k$ and a permutation $\sigma\in{\frak S}_n$,
$ka_{\sigma(i)} = a'_i$ ($1\le i \le n$) and $kh=h'$ hold.
\enddefinition

For each equivalence class of weight systems
there is exactly one reduced weight system satisfying
$a_1\le a_2\le \dots \le a_n$,
which we usually use in this paper.

\chap{2.3. Duality for weight systems.}

Let $n$ be a positive integer.
Take a weight system $W_a$ = $(a_1,\dots,a_n;h)$.
Assume that $h\in\sum_{i=1}^n \N a_i$ as always,
and that $W_a$ is reduced, for the sake of brevity.
We denote the integral vector ${}^t(a_1,\dots,a_n)$ by $a$.
\definition{Definition 2.3.1}
$a_0 := h - \sum_{i=1}^n a_i$.
\enddefinition

Assume also that $a_0 \ne 0$.

For a rational monomial $X_0^{m_0}\cdots X_n^{m_n}$ of degree $h$,
we assign an $n$-tuple of integers
$(\alpha_1, \dots, \alpha_n)$ = $(m_1 - 1, \dots, m_n - 1)$.
The whole such $n$-tuples constitute a set
\definition{Definition 2.3.2}
$$M(W_a) := \{(\alpha_1, \dots, \alpha_n) \in \Z^{\oplus n} \,|\, \sum_{i=1}^n
a_i(\alpha_i + 1) \equiv h \mod a_0\}.$$
\enddefinition

\proclaim{Lemma 2.3.3}
$M(W_a)$ is a subgroup of $\Z^{\oplus n}$ of index $|a_0|$.
\endproclaim
\demo{Proof}
This follows from
$M(W_a) = \{(\alpha_1,\dots, \alpha_n)\in\Z^{\oplus n}\,|\,\sum_{i=1}^n
a_i\alpha_i \equiv 0 \mod a_0\}$ and
$\gcd(a_1,\dots,a_n)=1$.
\qed
\enddemo

\remark{Remark 2.3.4}
If $W_a$ and $W_b$ are equivalent,
$M(W_a)$ and $M(W_b)$ are the same subset of $\Z^{\oplus n}$
up to the permutation of the coordinates.
We sometimes abbreviate this set as $M$ in this chapter.
\endremark
Let $C = (c_{ij})$ be an $n\times n$-matrix whose elements are nonnegative
integers.
Let $B$ be the $n\times n$-matrix $(c_{ij}-1)$.

\proclaim{Lemma 2.3.5}
Assume $Ca={}^t(h,\dots,h)$. Then
\roster
\item $(\det C)/h = (\det B)/a_0$, and this is an integer.
\item The following three conditions are equivalent:
(a) $\{(c_{i1}-1, \dots, c_{in}-1)|\, 1\le i \le n\}$ is a basis of $M$,
(b) $|\det B| = |a_0|$ and (c) $|\det C| = h$.
\endroster
\endproclaim
\demo{Proof}
Note that for $1\le i \le n$,
$(c_{i1}-1, \dots, c_{in}-1)\in \pi(a_0)\cap M$
and $(c_{i1}, \dots, c_{in})\in \pi(h)$,
where $\pi(t)$ is the hyperplane defined by $\sum_{i=1}^n a_i\alpha_i = t$.
The rest is clear.
\qed
\enddemo

\remark{Remark 2.3.6}
We will say a few words for the case $W_a$ is not reduced.
Let $d_0$ be $\gcd(a_1, \dots, a_n)$.
Then $M$ is a subgroup of index $|a_0|/d_0$, and the row vectors of $B$ is
a basis $\Longleftrightarrow |\det B| = |a_0|/d_0$
$\Longleftrightarrow |\det C| = h/d_0$.
\endremark

Let $W_a = (a_1, \dots, a_n;h)$ and $W_b = (b_1,\dots, b_n;k)$
be two weighted systems.

\definition{Definition 2.3.7}
An integer matrix $C \in M_n(\N)$ is said to be
a {\it weighted magic square\/} of weight $(W_a;W_b)$
if $C{}^t(a_1,\dots,a_n)={}^t(h,\dots,h)$ and
$(b_1,\dots,b_n)C = (k,\dots,k)$.
\enddefinition

\remark{Remark 2.3.8}
In the case $a_1 = \cdots = a_n = b_1 = \cdots = b_n = 1$ and $h = k$,
$C$ is called a {\it integer stochastic matrix\/} or {\it magic square\/},
which appears in classical combinatorics theory\cite{St}.
\endremark

\definition{Definition 2.3.9}
%$C$ is said to be {\it nondegenerate\/} if $\det C \ne0$.
A weighted magic square $C$ is
{\it primitive\/}
if $|\det C|$ = $h/\gcd(a_1,\dots,a_n)$ = $k/\gcd(b_1,\dots,b_n)$.
\enddefinition

\remark{Remark 2.3.10}
$C$ is primitive if and only if
the row vectors of $B$ span $M(W_a)$ and the column vectors span $M(W_b)$.
\endremark

\definition{Definition 2.3.11}
We say $W_a$ and $W_b$ are {\it dual\/} if
there exists a primitive weighted magic square $C$ of weight $(W_a;W_b)$.
Dual weight systems are {\it strongly dual\/}
if all rows and columns of $C$ contain $0$.
\enddefinition

\remark{Remark 2.3.12}
If $W_a$ and $W_b$ are reduced dual weights,
it follows that $h = k$,
$\sum_{i=1}^n a_i = \sum_{i=1}^n b_i$ and $a_0 = b_0$.

A permutation on weights $(a_1,\dots,a_n)$ (resp. $(b_1,\dots,b_n)$)
interchanges the corresponding columns (resp. rows) of $C$.

{}From this definition, one can calculate the dual weights.
\endremark

The following Proposition is handy for the calculation.

\proclaim{Proposition 2.3.13}
Let $W_a=(a_1,\dots,a_n;h)$ and $W_b=(b_1,\dots,b_n;h)$
be dual weight systems with $a_0$ be $h - \sum_{i=1}^n a_i > 0$
and $C = (c_{ij})$ be the corresponding weighted magic square.
Let $\Theta$ be the $(n-1)$-simplex in $M(W_a)_\Q$ with vertices
$\{(c_{i1}-1, \dots, c_{in}-1)|1\le i \le n\}$.
Then
\roster
\item $\sum_{j=1}^n c_{ij}a_j = h$ for $1 \le i \le n$,
\item $\displaystyle\frac{a_0}{h-a_0}(1,\dots,1) \in \Int\Theta$ and
\item $\Theta$ is elementary, i.e.,
$\Theta\cap M(W_a)$ are the set of all vertices of $\Theta$.
\endroster
Conversely, for a given $(a_1,\dots,a_n)\in(\Z)^{\oplus n}$,
if there exists
a $(n-1)$-simplex $\Theta$ with vertices
$\{(c_{i1}-1,\dots,c_{in}-1)|1\le i \le n\}$
which satisfies (1),(2) and (3) above,
then there exists a weight system $W_b = (b_1,\dots,b_n;h)$
such that $C$ is a weighted magic square of weight $(W_a;W_b)$.
If $W_b$ is reduced, then it is a dual weight system.
\endproclaim

\demo{Proof}
The assertions (1) and (3) are trivial.
The remaining (2) is equivalent to that
the origin sits inside of the $n$-simplex $\Delta$
which is the cone over $\Theta$ with apex $(-1,\dots,-1)$.
This is equivalent to that all the coefficients of the linear relation
between the vertices of $\Delta$ has a same sign.
Since we can take nothing but $b_i$'s as coefficients:
$b_0{}^t(-1,\dots,-1) + \sum_{i=1}^n b_i{}^t(c_{i1}-1,\dots,c_{in}-1) =0$,
(2) follows.

We will show the converse.
{}From the argument above,
$b_i$'s are determined as positive integers up to ratio.
We will take the reduced $b_i$'s at first.
We define $k := \sum_{i=0}^n b_i$ and $W_b :=(b_1,\dots,b_n;k)$.
Then $C$ is a weighted magic square of weight $(W_a;W_b)$.
We note that $|\det C|$ is a multiple of $k$, and by (3), $|\det C| = h$.
Thus multiplying some integer to $b_i$'s, we can achieve that $h = k$.
In particular, if this new $W_b$ is reduced
then $C$ is primitive since $k = |\det C|$.
\qed
\enddemo

\definition{Definition 2.3.14}
A weight system $W$ is {\it self-dual\/} if
$W$ and $W$ are dual.
\enddefinition

\remark{Remark 2.3.15}
If one can take a symmetric matrix as $C$, then $W$ is self-dual.
In general, there may be several dual weights as shown
in the Proposition below.
Similarly, self-duality does not mean that it is the only dual weight.
They are determined completely by a finite check.
\endremark

\chap{2.4. Examples.}

Here are the examples of case $a_0 = -1$,
which is hopefully to be generalized as a mirror of log Calabi-Yau manifolds.
\proclaim{Proposition 2.4.1}
$A_{l-1}=\{(1,k,l-k;l),\, 1 \le k \le l/2)\}$ $(l \ge 2)$ are
closed under duality for each $l$.
$D_{l+1}=(2,l-1,l;2l)$ ($l\ge3$),
$E_6=(3,4,6;12)$, $E_7=(4,6,9;18)$ and $E_8=(6,10,15;30)$ are self-dual.
\endproclaim
\demo{Proof}
Case $A_{l-1}$.

For an arbitrary positive integer $j$ such that $jk < l$, we can take $C$ as
$\left(\matrix l-jk & j & 0 \\ k & 0 & 1\\ 0 & 1 & 1\endmatrix\right)$.
Then the dual weight is $(1, j, l-j; l)$.
These are the only choice for $C$ up to permutation.

For $D_{l+1}$($l$:odd), $D_{l+1}$($l$:even), $E_6$, $E_7$ and $E_8$,
take the following matrices as $C$:
$$\left(\matrix (l+1)/2 & 1 & 0\\1 & 2 & 0\\ 0 & 0 & 2\endmatrix\right), \,
\left(\matrix l/2 & 0 & 1\\0 & 0 & 2\\ 1 & 2 & 0\endmatrix\right), \,
\left(\matrix 0 & 0 & 2\\0 & 3 & 0\\ 2 & 0 & 1\endmatrix\right), \,
\left(\matrix 0 & 3 & 0\\3 & 1 & 0\\ 0 & 0 & 2\endmatrix\right) \text{ and }
\left(\matrix 5 & 0 & 0\\0 & 3 & 0\\ 0 & 0 & 2\endmatrix\right).$$
They are all symmetric and unique up to the permutation of the rows.
They are strongly dual.
\qed
\enddemo

Next, we will treat the case of the exceptional unimodal
singularities\cite{AGV}.
They are all weighted homogeneous with $a_0 = 1$.

\proclaim{Theorem 2.4.2}
Let $W = (a_1,a_2,a_3;h)$ be a weight system which corresponds to
one of the 14 exceptional unimodal singularities.
Then $W$ has a unique dual weight system $W^*$ up to equivalence.
And $W^*$ is the weight system which corresponds to Arnold's strange duality.
\endproclaim

\demo{Proof}
We list below the choice of $C$.
These are, in fact, unique choice for primitive $C$ up to permutation.
They are automatically strongly dual.

Let us take coordinate as $x$, $y$ and $z$.
In the column of $C$, each monomial
$x^\alpha y^\beta z^\gamma$ means a row vector
$(\alpha, \beta, \gamma)$ of $C$.
For example, in the case $E_{14}$,
$C= \left(\matrix 4 & 0 & 1\\ 0 & 3 & 0\\ 0 & 0 &2\endmatrix\right)$,
which is nothing but the transposed matrix ${}^tC$ for $Q_{10}$.

We list also $C_0$,
which shows the vertices of the maximum Newton polytope of
polynomials of degree $h$ of $x$, $y$ and $z$,
and the Gabrielov numbers and the Dolgachev numbers\cite{AGV}.

\chap{Table of Arnold's singularities.}
\nopagebreak
$$\vbox{
\offinterlineskip
\halign{\strut\vrule\hfil#\vrule&\hfil{ #}&\hfil{ # }&\hfil{# }\vrule&\hfil{ #}
\vrule&\hfil{ #}&\hfil{ # }&\hfil{# }
\vrule&\hfil{ #}&\hfil{ # }&\hfil{# }
\vrule&&\hfil{ # }\vrule\cr
\noalign{\hrule}
class & $a_1$&$a_2$&$a_3$&$h$& &$C_0$&&&$C$&& $Gab$ &$Dol$ \cr
\noalign{\hrule}
$E_{12}$ & 6&14&21&42& $x^7$&$y^3$&$z^2$ & $x^7$&$y^3$&$z^2$ & 2 3 7& 2 3 7\cr
\noalign{\hrule}
$E_{13}$ & 4&10&15&30& $x^5y$&$y^3$&$z^2$ & $x^5y$&$y^3$&$z^2$ & 2 3 8& 2 4
5\cr
$Z_{11}$ & 6&8&15&30& $x^5$&$xy^3$&$z^2$ & $x^5$&$xy^3$&$z^2$ & 2 4 5& 2 3 8\cr
\noalign{\hrule}
$E_{14}$ & 3&8&12&24& $x^8$&$y^3$&$z^2$ & $x^4z$&$y^3$&$z^2$ & 2 3 9& 3 3 4\cr
$Q_{10}$ & 6&8&9&24& $x^4$&$y^3$&$xz^2$ & $x^4$&$y^3$&$xz^2$ & 3 3 4& 2 3 9\cr
\noalign{\hrule}
$Z_{12}$ & 4&6&11&22& $x^4y$&$xy^3$&$z^2$ & $x^4y$&$xy^3$&$z^2$ & 2 4 6& 2 4
6\cr
\noalign{\hrule}
$W_{12}$ & 4&5&10&20& $x^5$&$y^4$&$z^2$ & $x^5$&$z^2$&$y^2z$ & 2 5 5& 2 5 5\cr
\noalign{\hrule}
$Z_{13}$ & 3&5&9&18&  $x^6$&$xy^3$&$z^2$ & $x^3z$&$xy^3$&$z^2$ & 2 4 7& 3 3
5\cr
$Q_{11}$ & 4&6&7&18& $x^3y$&$y^3$&$xz^2$ & $x^3y$&$y^3$&$xz^2$ & 3 3 5& 2 4
7\cr
\noalign{\hrule}
$W_{13}$ & 3&4&8&16& $x^4y$&$y^4$&$z^2$ & $x^4y$&$z^2$&$y^2z$ & 2 5 6& 3 4 4\cr
$S_{11}$ & 4&5&6&16& $x^4$&$xz^2$&$y^2z$ & $x^4$&$xz^2$&$y^2z$ & 3 4 4& 2 5
6\cr
\noalign{\hrule}
$Q_{12}$ & 3&5&6&15& $x^5$&$y^3$&$xz^2$ & $x^3z$&$y^3$&$xz^2$ & 3 3 6& 3 3 6\cr
\noalign{\hrule}
$S_{12}$ & 3&4&5&13& $x^3y$&$xz^2$&$y^2z$ & $x^3y$&$xz^2$&$y^2z$ & 3 4 5& 3 4
5\cr
\noalign{\hrule}
$U_{12}$ & 3&4&4&12& $x^4$&$y^3$&$z^3$ & $x^4$&$y^2z$&$yz^2$ & 4 4 4& 4 4 4\cr
\noalign{\hrule}
}}$$
\qed
\enddemo

\remark{Remark 2.4.3}
K. Saito defines a duality of weights using a duality of
poset diagrams coming from the eigenvalues of the monodromy of
minimally elliptic singularities\cite{S4}.
His duality also reproduces Arnold's duality and he also computed
the dual weights for the 49 weight systems corresponding to
the minimally elliptic singularities which are not simple elliptic.

As shown below, our duality for $n = 3$ coincides with Saito's in the case
(0) the fourteen unimodal singularities and (3) $a_0 \ge 2$.
For cases (2) $a_0 = 1$ and modality $m$ is more than one,
our duality gives new dual weights.
These of course also correspond to the polar duality, as shown in Chapter 5.

It is natural to treat the simple elliptic cases (1) for $n=2$,
since these correspond to weighted elliptic curves.
Their list is the same as that of the weighted 2-spaces which correspond to
reflexive polytopes.
The weights are self-dual for $n=2$.

Notice that the singularities of dual weights are not dual itself
as Arnold's duality
but its (explicit) specializations enjoy the duality as shown in the later
Chapters.
\endremark

\proclaim{Proposition 2.4.4}
\roster
\item $(1,1;3)$, $(1,2;4)$, $(2,3;6)$ are self-dual.
\item $(2,2,3;8)$, $(2,2,5;10)$, $(2,3,4;10)$, $(2,4,7,14)$ and $(2,6,9,18)$
are self-strongly dual; $(2,3,6;12)$ and $(2,4,5;12)$ are strongly dual.
$(2,3,3;9)$ has no strongly dual weight.
\item Duality of weights for minimally elliptic weight systems with $a_0 > 1$
coincides with Saito's duality\cite{S4}.
\endroster
\endproclaim

\demo{Proof}
Here we give the full list of primitive $C$ up to permutation.
The notation of classes is after \cite{AGV} and \cite{HYW}.

\chap{(1) Simple elliptic singularities.}
{\eightpoint
$\vbox{
\offinterlineskip
\halign{\strut\vrule#\vrule& \hfil # \hfil&\hfil # \hfil\vrule
&& \hfil#\hfil\vrule\cr
\noalign{\hrule}
cl. & $C$& & $a_1$ $a_2$ & $h$ & dual \cr
\noalign{\hrule}
$P_8$ & $x^2y$ & $xy^2$ & 1 1 & 3 & 1 1 \cr
$X_9$ & $y^2$ & $x^2y$ & 1 2 & 4 & 1 2 \cr
$J_{10}$ & $x^3$ & $y^2$ & 2 3 & 6 & 2 3 \cr
\noalign{\hrule}
}}$
}

\chap{(2) Minimal elliptic singularities with $a_0=1$ and $m>1$.}
{\eightpoint
$\vbox{
\offinterlineskip
\halign{\strut\vrule#\vrule& \hfil # &\hfil # \hfil&\hfil # \hfil\vrule
&& \hfil#\hfil\vrule\cr
\noalign{\hrule}
cl. & &$C$& & $a_1$ $a_2$ $a_3$ & $h$ & dual \cr
\noalign{\hrule}
$V_{1,0}=V_{15}$ & $xz^2$ & $yz^2$ & $x^2y^2$ & 2 2 3 & 8 & 2 2 3 \cr
$              $ & $x^2y^2$ & $xy^3$ & $xz^2$ &       & * & 1 2 4 \cr
\hphantom{$N_{1,0}=$M}$N_{16}$ & $x^3y^2$ & $x^2y^3$ & $z^2$ & 2 2 5 & 10 & 2 2
5 \cr
$U_{1,0}=U_{14}$ & &-& & 2 3 3 & 9 & - \cr
$              $ & $y^2z$ & $x^3y$ & $yz^2$ &       & * & 1 3 4 \cr
$S_{1,0}=S_{14}$ & $y^2z$ & $x^2y^2$ & $xz^2$ & 2 3 4 & 10 & 2 3 4 \cr
$              $ & $xz^2$ & $x^3z$ & $y^2z$ &       & * & 1 3 5 \cr
$W_{1,0}=W_{15}$ & $y^2z$ & $x^3y^2$ & $z^2$ & 2 3 6 & 12 & 2 4 5 \cr
$              $ & $z^2$ & $x^3z$ & $y^2z$ &       & * & 1 4 6 \cr
$Q_{2,0}=Q_{14}$ & $y^3$ & $x^2y^2$ & $xz^2$ & 2 4 5 & 12 & 2 3 6 \cr
$Z_{1,0}=Z_{15}$ & $xy^3$ & $x^3y^2$ & $z^2$ & 2 4 7 & 14 & 2 4 7 \cr
$J_{3,0}=J_{16}$ & $y^3$ & $x^3y^2$ & $z^2$ & 2 6 9 & 18 & 2 6 9 \cr
\noalign{\hrule}
}}$
}

The * signifies that the dual weight is not strongly dual,
hence the singularity corresponding to the smaller triangle is reducible.
Other weight systems are strongly dual.

\chap{(3) Minimal elliptic singularities with $a_0>1$.}
{\eightpoint
$\vbox{
\offinterlineskip
\halign{\strut\vrule#\vrule& \hfil # &\hfil # \hfil&\hfil # \hfil\vrule
&& \hfil#\hfil\vrule\cr
\noalign{\hrule}
cl. & &$C$& & $a_1$ $a_2$ $a_3$ & $h$ & dual \cr
\noalign{\hrule}
$V'_{18}$ & &-& & 3 3 4 & 12 & - \cr
$U_{16}$ & $x^5$ & $y^2z$ & $yz^2$ & 3  5  5 & 15 & 3  5  5 \cr
$S_{16}$ & $x^4y$ & $xz^2$ & $y^2z$ & 3  5  7 & 17 & 3  5  7 \cr
$W_{17}$ & $x^5y$ & $z^2$ & $y^2z$ & 3  5 10 & 20 & - - - \cr
$Q_{16}$ & $x^4z$ & $y^3$ & $xz^2$ & 3  7  9 & 21 & 3  7  9 \cr
$Z_{17}$ & $x^4z$ & $xy^3$ & $z^2$ & 3  7 12 & 24 & - - - \cr
$E_{18}$ & $x^5z$ & $y^3$ & $z^2$ & 3 10 15 & 30 & - - - \cr
\noalign{\hrule}
}}$
}
{\eightpoint
$\vbox{
\offinterlineskip
\halign{\strut\vrule#\vrule& \hfil # &\hfil # \hfil&\hfil # \hfil\vrule
&& \hfil#\hfil\vrule\cr
\noalign{\hrule}
${}_2V^*_{18}$ & $xy^3$ & $x^3z$ & $yz^2$ & 4  5  7 & 19 & 4  5  7 \cr
${}_1V^*_{18}$ & $x^3z$ & $y^4$ & $xz^2$ & 4  5  8 & 20 & 4  5  8 \cr
$N_{19}$ & $x^3z$ & $xy^4$ & $z^2$ & 4  5 12 & 24 & - - - \cr
$S_{17}$ & $x^6$ & $xz^2$ & $y^2z$ & 4  7 10 & 24 & - - - \cr
$W_{18}$ & $x^7$ & $y^2z$ & $z^2$ & 4  7 14 & 28 & 4  7 14 \cr
$Q_{17}$ & $x^5y$ & $y^3$ & $xz^2$ & 4 10 13 & 30 & - - - \cr
$Z_{18}$ & $x^6y$ & $xy^3$ & $z^2$ & 4 10 17 & 34 & 4 10 17 \cr
$E_{19}$ & $x^7y$ & $y^3$ & $z^2$ & 4 14 21 & 42 & - - - \cr
\noalign{\hrule}
}}$
}

{\eightpoint
$\vbox{
\offinterlineskip
\halign{\strut\vrule#\vrule& \hfil # &\hfil # \hfil&\hfil # \hfil\vrule
&& \hfil#\hfil\vrule\cr
\noalign{\hrule}
${}_3V^*_{19}$ & $y^4$ & $x^3z$ & $yz^2$ & 5  6  9 & 24 & - - - \cr
${}_2N_{20}$ & $y^5$ & $x^3z$ & $z^2$ & 5  6 15 & 30 & - - -\cr
\noalign{\hrule}
}}$
}

{\eightpoint
$\vbox{
\offinterlineskip
\halign{\strut\vrule#\vrule& \hfil # &\hfil # \hfil&\hfil # \hfil\vrule
&& \hfil#\hfil\vrule\cr
\noalign{\hrule}
$V'_{20}$ & $x^3z$ & $z^3$ & $xy^3$ & 6  7  9 & 27 & 6  7  9 \cr
${}_2V^*_{19}$ & $x^5$ & $xy^3$ & $yz^2$ & 6  8 11 & 30 & - - - \cr
${}_1V^*_{19}$ & $x^4y$ & $y^4$ & $xz^2$ & 6  8 13 & 32 & 4  7 16 \cr
${}_1N_{20}$ & $x^5y$ & $xy^4$ & $z^2$ & 6  8 19 & 38 & 6  8 19 \cr
$Q_{18}$ & $x^8$ & $y^3$ & $xz^2$ & 6 16 21 & 48 & 3 16 24 \cr
$Z_{19}$ & $x^9$ & $xy^3$ & $z^2$ & 6 16 27 & 54 & 4 18 27 \cr
$E_{20}$ & $x^{11}$ & $y^3$ & $z^2$ & 6 22 33 & 66 & 6 22 33 \cr
\noalign{\hrule}
}}$
}
{\eightpoint
$\vbox{
\offinterlineskip
\halign{\strut\vrule#\vrule& \hfil # &\hfil # \hfil&\hfil # \hfil\vrule
&& \hfil#\hfil\vrule\cr
\noalign{\hrule}
$V'_{21}$ & $z^3$ & $y^4$ & $x^3z$ & 8  9 12 & 36 & 8  9 12 \cr
$V^*_{20}$ & $y^4$ & $x^5$ & $yz^2$ & 8 10 15 & 40 & 5  8 20 \cr
$N_{21}$ & $y^5$ & $x^5y$ & $z^2$ & 8 10 25 & 50 & 8 10 25 \cr
\noalign{\hrule}
}}$
}
\qed
\enddemo

\remark{Remark 2.4.5}
The four weight systems in (3):
$(6,8,13;32)$, $(6,16,21;48)$, $(6,16,27;54)$ and $(8,10,15;40)$
are not in the list of Reid\cite{Re}\cite{F},
hence do not correspond to K3 surface with cyclic quotient singularities.

For the weights with some $C$ but without a dual weight (e.g. $W_{17}$),
$C$ does not satisfy primitivity for the dual weight.
On the contrary, they enjoy a dual correspondence defined in Chapter 1,
except $(5,6,15;30)$.

All the dual weights are strongly dual.
\endremark

\remark{Remark 2.4.6}
There are 41 weight systems with $a_0=1$ in Reid's list and
19 of them are not minimally elliptic.
A similar calculation shows that all the duality for them are as follows;
$(1,1,1;4)$, $(1,1,2;5)$, $(1,2,2;6)$, $(1,2,3;7)$ and $(1,2,4;8)$ are
self-dual;
$(1,1,3;6)$ and $(1,2,2;6)$ are dual;
$(1,2,4;8)$, $(1,3,4;9)$, $(1,3,5;10)$ and $(1,4,6;12)$ have dual weight
systems which are minimally elliptic as in Proposition 2.4.4 (2);
and all the ten weight systems with $a_1>1$ do not have dual weight systems.
\endremark

\chap{3. Relation of duality of weight systems and polar duality.}

We shall show in this chapter that
the duality of weight for $a_0 = 1$ is a special case of
dual correspondence defined in Chapter 1.
We will always choose one coordinate for
the compactification parameter.

Let $W_a = (a_1,\dots,a_n;h)$ be a weight system with
$h \in \sum_{i=1}^n \N a_i$
and $a_0 := h - \sum_{i=1}^n a_i$ as usual.
We assume $a_0 > 0$ in this section.

Let $\P(a)$ be the weighted projective space of weight $(a_0,a_1,\dots,a_n)$
and $X_0, \dots, X_n$ be the natural homogeneous coordinates of $\P(a)$
induced by the definition.
Then the $n$-dimensional complex torus $\T$ is embedded equivariantly
in $\P(a)$ as the locus $\{F_0 := \prod_{i=0}^n X_i\ne0\}$.

Let $M(W_a)$ be the group defined in Section 3.
$M(W_a)$ is nothing but the abelian group of exponents of rational monomials,
which appears in the context of toric construction.
For a point $P=(\alpha_1,\dots,\alpha_n)$ on $M(W_a)$,
we write the corresponding monic monomial by
$F_P=\prod_{i=0}^nX_i^{\alpha_i+1}$,
where $\alpha_0$ is determined by
$\sum_{i=0}^n a_i\alpha_i=0$.

\definition{Definition 3.1}
Let $\overline\Delta(a)$ be the $n$-simplex in $M(W_a)_\Q$
which is the convex hull of the following vertices:
$(-1+h/a_1,-1,\dots,-1)$, $\ldots$,
$(-1,\dots,-1,-1+h/a_n)$ and $(-1,\dots,-1)$.
\enddefinition

$(\overline\Delta, \overline M)$ =
$(\overline\Delta(a), M(W_a))$ is an element of ${\Cal P}_n$
with $\P_{\overline\Delta,\overline M} \cong \P(a)$.

\proclaim{Lemma 3.2}
Let $(\overline\Delta,M)$ be as above.
Then its dual $(\overline\Delta^*, M^*)$ is integral
and the vertices generate $M^*$.
\endproclaim

\demo{Proof}
For the sake of simplicity, we assume that $W_a$ is reduced.
Recall that $M$ =
$\{(\alpha_1,\dots,\alpha_n)|\sum_{i=1}^n a_i\alpha_i \equiv 0 \mod a_0\}$.
Let $M' := \Z^{\oplus n}$, which contains $M$
as a $\Z$-submodule of index $a_0$.
If one computes the dual of $(\overline\Delta, M')$,
it is an $n$-simplex with vertices $v_1 := (1,0,\dots,0)$, $\ldots$,
$v_n:=(0,\dots,0,1)$ and $v_0 := (-a_1/a_0,\dots,-a_n/a_0)$ in $(M')^*_\Q$.
Thus $v_1$, $\ldots$, $v_n$ generates $(M')^*$.
Since $(M')^*$ is a submodule of $M^*$,
Those $n$ vertices are also contained in $M^*$.
On the other hand, $v_0$ also belongs to $M^*$
since $\sum_{i=1}^n(-a_i/a_0)\alpha_i \in \Z$ for any
$(\alpha_1,\dots,\alpha_n) \in M$.

The order of $v_0$ in $M^*/(M')^*$ is exactly $a_0$ since
$(a_1,\dots,a_n)$ is reduced.
Thus the last statement follows.
\qed
\enddemo

\remark{Remark 3.3}
By the previous proof, the submodule $(M')^*$ generated by
$v_1,\dots,v_n$ has index $a_0$ in $M$.
\endremark

Let $(\overline\nabla,\overline N)$ =
$(\overline\Delta(b), M(W_b))$ be another pair.
We assume $W_b$ is also reduced for the sake of simplicity.

We restrict here the correspondence $\sigma$ in Chapter 1 as follows:
(*)
\roster
\item $\sigma$ sends
the vertex of $(\overline\Delta)^*$ which corresponds to
the divisor $\{X_0=0\}$
to the vertex of $\overline\nabla$ which corresponds to the point
$(1:0:\cdots:0)$ (i.e., $X_1=\cdots=X_n=0$).
\item same for ${}^t\sigma$.
\endroster

\proclaim{Proposition 3.4}
Assume that $(\overline\Delta,\overline M)$ and
$(\overline\nabla,\overline N)$ has a dual correspondence $\sigma$
which satisfies the condition (*) above.

Let $(c_{i1}-1,\dots,c_{in}-1)$ ($1\le i \le n$)
be the vertices except $(-1,\dots,-1)$ of $\sigma({\overline\nabla}^*)$
in $\overline M$.
Let $B$ be the $n\times n$-matrix $(c_{ij}-1)$.

Then
\roster
\item $a_0|h$, $b_0|h$ and $|\det B|=a_0b_0$,
\item If $a_0 = b_0 =1$, $W_a$ and $W_b$ are dual weight systems.
\endroster
\endproclaim

\demo{Proof}
(1) By $\sigma({\overline N}^*) = \overline M$ and
$\sigma(-b_1/b_0,\dots,-b_n/b_0)=(-1,\dots,-1)$,
$(-1,\dots,-1) \in \overline M$. Thus $a_0|h$.
Similarly for $b_0|h$.
Since $\sigma(1,0,\dots,0)$, $\ldots$, $\sigma(0,\dots,0,1)$ generates
a submodule of index $b_0$ in $\overline M$,
$|\det B|=a_0b_0$.

(2) Let $C$ be $(c_{ij})$.
By condition (*), $C$ is a weighted magic square of weight $(W_a;W_b)$ and
primitive by $|\det B| = 1$ using Lemma 2.3.5.
\qed
\enddemo

Thus we have

\proclaim{Theorem 3.5}
Assume two weight systems $W_a$ and $W_b$ with $a_0 = b_0 =1$ are dual.
Then
\roster
\item
there exists an identification $M = (\overline N)^* = \overline M$ such that
${\overline\nabla}^*$ is an integral simplex and
a subcone of $\overline\Delta$ with a common apex $(-1,\dots,-1)$ and
a common plane $\pi$ = $\{\sum_{i=1}^na_i\alpha_i = a_0\}$ of facet.
The vertices of ${\overline\nabla}^*$  on $\pi$ generates $M$.

Similarly for the induced identification $N = (\overline M)^* = \overline N$.
\item If a linear subsystem of anticanonical divisor of $\P(a)$ corresponds to
some reflexive polytope $\Delta$ such that
${\overline\nabla}^* \subset \Delta \subset \overline\Delta$,
then the family corresponds to $\Delta^*$ is a linear subsystem of
anticanonical divisor of $\P(b)$.
This relation is inclusion-reversing with respect to polytopes.
\endroster
\endproclaim

\remark{Remark 3.6}
The facet above is the Newton polytope for
an $(n-1)$-dimensional hypersurface singularity $X_0$ with $\C^*$-action.
$\bar Z$ is the toric compactification in $\P(a)$ of a deformation $Z$ of
$X_0$.
\endremark

\remark{Remark 3.7}
Let $\Delta$ be the {\it full Newton polytope\/} of $\P(a)$,
that is the convex hull of all vertices corresponding to
the monomial anticanonical divisors.
We have the associated projective variety $\P_\Delta$ as usual.
In general, $\Delta$ does not contain all the generators of
the ring of regular functions of $\P(a)$.
Thus $\P(a)$ and $\P_\Delta$ are not isomorphic in general.
\endremark

\remark{Remark 3.8}
By \cite{B1, Theorem 5.4.5, Corollary 5.4.6},
for a reflexive simplex $\Delta = \overline\Delta$,
that is, for the case $a_i|h$ for all $i$,
$\P_\Delta$ is a weighted projective variety whose weight is the
coefficients
$(a_0,\dots,a_n)$ of the unique linear relation
$\sum_{i=0}^n a_iP_i = 0$ ($a_i\in\Z_+$, $\gcd(a_0,\dots,a_n) = 1$)
among the vertices $\{P_i\}$.
Moreover, the polar dual family of deformation of Fermat-type hypersurfaces
is the quotient of a subfamily by $\pi_1(\Delta,M)$ \cite{B1 Corollary 5.5.6}.
\endremark

\remark{Remark 3.9}
Borcea's construction\cite{Bo} does not assume the property (*) which is
coming from singularity theory.
As he has pointed out, for a given reflexive polytope,
its realization as a family of weighted hypersurface has some choice.
Nevertheless, for Arnold's singularities,
our duality of weights is unique and compatible with
that of usual singularity theory.
\endremark

\chap{4. Application for $n=3$.}
\chap{4.1. K3 surfaces.}
We refer the reader \cite{BPV} for more complete reference about K3 surfaces
and \cite{AM} for general arguments about string theory on K3 surfaces.

Let $X$ be K3 surface, that is, a compact complex surface
with trivial canonical bundle and irregularity $h^1(\CO_X) = 0$.
It is well known that any K3 surface is
a simply connected, K\"ahler and symplectic manifold.
The only nonzero Hodge numbers of $X$ are
$h^{0,0}=h^{2,0}=h^{0,2}=h^{2,2}=1$ and $h^{1,1}=20$.

First we recall here some notions about lattices.
\definition{Definition 4.1.1}
Let $L$ be a free $\Z$-module of finite rank
with a symmetric bilinear $\Z$-valued pairing $\left<\, ,\, \right>$.
Such an $L$ is called a {\it lattice\/}.
$L$ is {\it even\/} if for each $x \in L$,
$\left<x,x\right>$ is an even integer,
and {\it unimodular\/} if the discriminant is $\pm1$.
The {\it index\/} of $L$ is the triplet $(\lambda_+,\lambda_0,\lambda_-)$
of the numbers of positive, null and negative eigenvalues
of the corresponding matrix.
$L$ is said to be {\it positive definite\/}, {\it negative definite\/},
{\it definite\/}, {\it indefinite\/}, respectively,
if the corresponding matrix is so.
\enddefinition

Since $\pi_1(X)=\{1\}$,
we can identify $H^2(X,\Z)$ and its natural image in $H^2(X,\C)$.
$H^2(X,\Z)$ is naturally regarded as a lattice by the cup product,
and it is an even unimodular lattice of index $(3,0,-19)$,
so is isomorphic to $L := (-E_8)^{\oplus2}\oplus U^{\oplus3}$.
Here, $-E_8$ is the even unimodular negative definite lattice of rank $8$,
whose matrix is $(-1)$ times
that of the lattice corresponding to the Dynkin diagram $E_8$.
$U$ is the even unimodular indefinite lattice of rank $2$,
whose matrix is $\left(\matrix0&1\\1&0\endmatrix\right)$.

In general,
the algebraic lattice $S_X$ and the transcendental lattice $T_X$ of
a compact complex surface $X$ are defined as sublattices
in (the image of) $H^2(X,\Z)$ as
$S_X := H^{1,1}(X,\R)\cap H^2(X,\Z)$ and
$T_X := S_X^\bot$ in $H^2(X,\Z)$.
Here $\bot$ denotes the orthogonal complement.
The rank of $S_X$ is called the {\it Picard number\/} of $X$
and denoted by $\rho(X)$.

%We call a pair of K3 surface and its K\"ahler class $(X,\kappa)$
%a {\it marked K3 surface\/}.

Ruan and Tian have constructed the {\it quantum cohomology ring\/}
for $\Z$-coefficient of symplectic manifolds in \cite{RT}.
As a graded $\Z$-module, the quantum cohomology ring is
same as the classical cohomology ring,
but in general, the multiplicative structures are different.
The mirror map $\mu$ is usually thought to be an involution of
moduli of Calabi-Yau manifolds with K\"ahler structure
which induces a kind of isomorphism on
the quantum cohomology rings of the Calabi-Yau manifolds exchanged by $\mu$.

The following is an easy modification of \cite{RT Example 8.4}.

\proclaim{Proposition 4.1.2}
Let $X$ be a compact K\"ahler surface with trivial canonical bundle.
Then the quantum cohomology ring of $X$ is isomorphic to
the classical cohomology ring of $X$.
\endproclaim

Thus for a K3 surface or a 2-dimensional complex torus,
we may assume that a mirror map preserves the usual cup product.

\chap{4.2. Fixed set of the mirror map.}

First we shall proof some finiteness result of a special automorphism of
marked K3 surface.

\proclaim{Theorem 4.2.1}
Let $X$ be a K3 surface with a K\"ahler form $\omega$,
$L$ be the second integral cohomology lattice $H^2(X,\Z)$ and
$G$ be a subgroup of $\Aut X$.
Assume that the induced action of $G$ on $L$
fixes a nonzero holomorphic $2$-form and $\omega$.
Then
\roster
\item The fixed lattice $L^G$ contains $T_X$.
\item The orthogonal complement of $(L^G)^\bot$ in
$L$ is a negative definite sublattice of $S_X$.
\item $G$ is a finite group.
\endroster
\endproclaim

\remark{Remark 4.2.2}
If the algebraic dimension $a(X)$ of $X$ is not one,
it is essentially proved in \cite{Ni}.
\endremark

First we note the following elementary fact.
\proclaim{Lemma 4.2.3}
Let $L$ be a lattice and $K$ be a subset of $L$.
Let $\phi$ be an isometry on $L$.
\roster
\item If $\phi(K)=K$ then $\phi(K^\bot) = K^\bot$.
\item If $\phi$ acts as identity on $K$, then for any $x \in L$,
$\phi(x)-x \in K^\bot$.
\endroster
\endproclaim

\proclaim{Corollary 4.2.4}
Let $X$ be a compact K\"ahler surface and $g \in \Aut X$.
Then $g^*(S_X) = S_X$ and $g^*(T_X) = T_X$.
\endproclaim

\demo{Proof of Theorem}
First note that the kernel of the multiplication of a nonzero holomorphic
two form $\Omega$ on $H^2(X,\R)$ is $H^{1,1}\cap H^2(X,\R)$.

Take an arbitrary element $x \in T_X$ and $g \in G$ .
Then $y := g^*x-x \in T_X$.
By the previous lemma,
$y \in T_X \cap \Omega^\bot \cap \omega^\bot$ =
$T_X \cap S_X \cap \omega^\bot$.
Since $y \in T_X \cap S_X$, we have $\left<y,y\right>=0$.
On the other hand, since $S_X \cap \omega^\bot$ is negative definite
by the signature theorem, one has $y = 0$.
Thus $g^*x =x$.
This proves $(1)$.

$(2)$ follows from
$(L^G)^\bot \subset (T_X)^\bot \cap \Omega^\bot \cap \omega^\bot$
= $S_X \cap \omega^\bot$.

Next we shall prove $(3)$.

Since $(L^G)^\bot$ is negative definite and finitely generated,
the isometry group $O((L^G)^\bot)$ is a finite group.

Let $j : G \to O((L^G)^\bot)$ be the natural group homomorphism.
We claim that $j$ is injective.

Suppose the induced action of an element $g \in G$ on $(L^G)^\bot$
is identity.
For an arbitrary $x \in L$, $y := g^*x-x$ belongs to $(L^G)^\bot$.
Thus $g^*y = y$.
Then for an integer $n$,
we have $(g^n)^*x = x+ ny$ and since $g^*$ is isometry,
$\left<(g^n)^*x, (g^n)^*x\right> = \left<x,x\right>$.
This leads to $y = 0$ by definiteness of $(L^G)^\bot$.
Thus $g^*$ is the identity on $L$.

Since $g^*$ is, by definition, an effective Hodge isometry on $L$,
$g$ is the identity element in $\Aut X$ by Torelli's theorem.
\qed
\enddemo

\remark{Remark 4.2.5}
The structure of $G$ can be found in the list of \cite{Ni}
once one knows that $G$ is a finite algebraic automorphism group.
Also, $(L^G)^\bot$ does not contain square $(-2)$-elements.
The lattice $L_0$ in the next section is nothing but $(L^G)^\bot$
when the K\"ahler class is generic in $(L_D)_\R$.
\endremark

\proclaim{Corollary 4.2.6}
$L^G + (L^G)^\bot$ is a sublattice in $L$ of finite index and
$L_\Q = (L^G)_\Q \oplus (L^G)^\bot_\Q$.
\endproclaim

\demo{Proof}
For any $x \in L$, there is an orthogonal decomposition in $\Z[1/|G|]$:
$$x = \frac{1}{|G|}(x_1 \oplus x_2), \,
x_1 = \sum_{g\in G}g^*x \in (L^G), \,
x_2 = \sum_{g\in G}(x-g^*x) \in (L^G)^\bot.\qed$$
\enddemo

\chap{4.3. Monomial divisor mirror map for K3 toric hypersurfaces.}

In this section, we will treat the algebraic K3 surfaces which are
the minimal resolutions of toric hypersurfaces.

\definition{Definition 4.3.1}
A {\it polyhedron\/} is a polytope of dimension three.
\enddefinition

For a given reflexive polyhedron $\Delta$,
we use notations $\Theta$ and $\Gamma$
for $2$-and $1$-dimensional face of $\Delta$.
The dual face $\Gamma^*$ of $\Gamma$ is $1$-dimensional face of $\Delta^*$,
so a summation over $\Gamma$ is same as that over $\Gamma^*$.

We denote by $Z$ the intersection of a $\Delta$-regular anticanonical section
of $\P_\Delta$ and the $n$-dimensional orbit of the torus,
and by $\bar Z$ its closure, namely the section itself.
$\bar Z$ has only canonical Gorenstein singularities
by Theorem 1.1.8 due to Batyrev.
Let $\tilde Z$ be the minimal resolution of $\bar Z$.
$X = \tilde Z$ is an algebraic K3 surface.

We define the following sublattices of $H^2(X,\Z)$ for
toric hypersurface K3 surface $X$.

\definition{Definition 4.3.2}
\roster
\item $L_G := H^{1,1}(H^2_c(Z,\C)) \cap H^2(X,\Z)$,
\item $L_D := Im(H^{1,1}(\tilde\P_\Delta) \to H^{1,1}(X)) \cap H^2(X,\Z)$,
\item $L_0 := L_D^\bot$ in $S_X$,
%$L_0 := \omega^\bot in S_X$ and
%$L_D := L_0^\bot in S_X$.
\endroster
where $\tilde\P_\Delta$ is a toric resolution of $\P_\Delta$.
%We call them {\it Gabrielov lattice\/}, {\it Dolgachev lattice\/} and
%{\it inner lattice\/}, respectively.
\enddefinition

Using a result of \cite{DK} or \cite{B2}, we get

\proclaim{Theorem 4.3.3}
Let $\Delta$ be three-dimensional reflexive polytope
and $\tilde Z$ be the K3 surface which is the minimal resolution of
a $\Delta$-regular anticanonical section of $\P_\Delta$. Then
\roster
%\item $L_G(\tilde Z) = (\Delta^{[1]} \cap M)/M
%\item $L_D(\tilde Z) = (\Delta^{*[1]} \cap N)/N
\item $\rk L_G(\tilde Z) = l(\Delta^{[1]}) - 3$,
\item $\rk L_D(\tilde Z) = l(\Delta^{*[1]}) - 3$,
\item $\rk L_0(\tilde Z) = \sum_{\dim\Gamma=1} l^*(\Gamma^*)l^*(\Gamma)$.
\endroster
\endproclaim

\demo{Proof}
By definition we can follow the same argument of calculation of
$H^{1,1}_{toric}$ in \cite{AGM},
one can show that
$L_D(\tilde Z) = (\Delta^{*[1]} \cap N)/N$ and in particular,
$\rk L_D(\tilde Z) = l(\Delta^{*[1]}) - 3$.

By \cite{DK 5.11},
$h^{1,1}(Z) = l^*(2\Delta)-4l^*(\Delta)-\sum_\Theta l^*(\Theta)$
= $l(\Delta) - 4 - \sum_\Theta l^*(\Theta)$
= $l(\Delta^{[1]}) - 3$.
Since the toric compactification $\T \to \P_\Delta$ adds the divisors
corresponding to the vertices of $\Delta^*$,

$h^{1,1}(\bar Z) = h^{1,1}(Z) + l(\Delta^{*[0]}) - 3$.

A toric resolution of the ambient space produces the divisors corresponding
to the integral points $P$ of $\Delta^* - (\Delta^{*(0)}\cup\{0\})$.
We denote by $C_P$ the support of the center of blowing up in $\P_\Delta$
corresponding to $P$.
If $P$ is in some interior of a $2$-dimensional face of $\Delta^*$,
$C_P$ is a point which is disjoint from $\bar Z$,
since $\bar Z$ is $\Delta$-regular.

If $P$ is in some $\Int\Gamma^*$,
$C_P$ is a curve in $\P_\Delta$
which intersects normally with $\bar Z$ at ($l^*(\Gamma) + 1$) points.
Since each toric blowing-up corresponding to $P$ creates one ambient divisor,
its orthogonal complement increases $L_0$ by $l^*(\Gamma)$.

Thus we have
$h^{1,1}(\tilde Z)$
= $h^{1,1}(\bar Z) + \sum_{\Gamma^*} (l^*(\Gamma)+1)l^*(\Gamma^*)$
= $l(\Delta^{[1]}) + l(\Delta^{*[1]}) - 6
+ \sum_{\Gamma^*} l^*(\Gamma)l^*(\Gamma^*)$
and
$\rk L_0(\tilde Z)$ = $\sum_{\Gamma} l^*(\Gamma)l^*(\Gamma^*)$.

We get
$h^{1,1}(\tilde Z) - \rk L_D - \rk L_0$
= $l(\Delta^{[1]}) - 3$
= $\rk (\Delta^{[1]} \cap M)/M$,
which is the expected dimension for $H^{1,1}_{poly}$ in \cite{AGM},
and this is nothing but $\rk L_G$ = $h^{1,1}(H^2_c(Z))$.
\qed
\enddemo

\proclaim{Corollary 4.3.4}
For a fixed reflexive polyhedron $\Delta$,
$L_G(\tilde Z)$, $L_D(\tilde Z)$ and $L_0(\tilde Z)$
are independent of the choice of $\Delta$-regular $Z$.
\endproclaim

We will thus write them as
$L_G(\Delta)$, $L_D(\Delta)$ and $L_0(\Delta)$, respectively.

\proclaim{Corollary 4.3.5}
\roster
\item $L_G(\Delta) = L_D(\Delta^*)$, $L_0(\Delta) = L_0(\Delta^*)$.
\endroster
\endproclaim

\remark{Remark 4.3.6}
Since the rank of $H^{1,1}$ is $20$,
$\rk L_G(\Delta) + \rk L_D(\Delta) + \rk L_0(\Delta) = 20$.
That formula gives a relation
$\sum_\Gamma (l^*(\Gamma)+1)(l^*(\Gamma^*)+1) = 24$.
\endremark

\definition{Definition 4.3.7}
We say the natural isomorphism $\mu\: (L_D(\Delta))_\C \to (L_G(\Delta^*))_\C$
a {\it monomial divisor mirror map\/} according to \cite{AGM}.
\enddefinition

\remark{Remark 4.3.8}
In fact, we can easily determine the lattices $L_D$ and $L_0$
from the data $(\Delta,M)$.
For instance, the dual graph of $L_D$ is obtained from $\Delta^{*[1]}$ by
removing three integral points which consist a basis of $N$.
Each remaining vertex is a $(-2)$-curve
and each remaining integral point on a edge $\Gamma$ is a sum of
$(l^*(\Gamma^*)+1)$ disjoint $(-2)$-curves.
Each edge joining integral points shows the corresponding cycles
intersect with an apparent intersection number.

We conjecture that the same rule of intersection product
holds for $L_G$ and $\Delta^{[1]}$.

In the case of Arnold's strange duality,
the lattice structure of $L_G(\Delta)$ is isomorphic to that of
$L_D(\Delta^*)$.

In general, $L_G\oplus U$ is not isomorphic as a lattice over $\Z$
to the Milnor lattice but it seems so over $\Z[1/G]$.
This may be for a future investigation.
\endremark

\proclaim{Theorem 4.3.9}
Let $X_0$ be one of Arnold's singularities
and $W_a = (a_1,a_2,a_3;h)$ be the corresponding weight system.
\roster
\item There exists a unique dual weight system $W_b$.
\item $W_b$ coincides with the weight system of the dual singularity
in the sense of Arnold.
\item The full Newton polyhedron, namely the convex hull of
$\overline\Delta(a) \cap M(W_a)$ in $M(W_a)_\Q$ is a reflexive polyhedron
containing $(\overline\nabla)^*$ with the dual correspondence above.
\item Let $\Delta$ be any reflexive polyhedron in $\overline\Delta(a)$
which contains $(\overline\nabla)^*$.
The intersection with $\Delta$ and the plane
$\sum_{i=1}^3 a_i\alpha_i=1$ is a facet which determines
the equivalent singularity to $X_0$.
\item $\Delta^*$ is a reflexive polyhedron in $\overline\nabla$ containing
$(\overline\Delta)^*$.
\item There exists a reflexive $\Delta$ such that
$L_0 = 0$ and $\mu$ interchanges $L_G(\Delta)$ and $L_D(\Delta) = S_X$.
\endroster
\endproclaim

\demo{Proof}
(1)(2)(5) is done in previous chapters.
(3) is well-known and (4) is easy case by case check.
Here we show an example of choices of $\Delta$ with $L_0=0$.
For dual singularities, these reflexive polytopes are dual to each other.
For self-dual singularities, these examples satisfy $\Delta^* \cong \Delta$.
These choices of $\Delta$ are not unique in general.
$$\vbox{
\offinterlineskip
\halign{\strut\vrule\hfil#\vrule&\hfil{ #}&\hfil{ # }&\hfil{# }\vrule&\hfil{ #}
\vrule&\hfil{ # }\hfil\vrule\cr
\noalign{\hrule}
class & $a_1$&$a_2$&$a_3$&$h$& $\Delta$\cr
\noalign{\hrule}
$E_{12}$ & 6&14&21&42& $W^{42}$, $X^7$, $Y^3$, $Z^2$\cr
\noalign{\hrule}
$E_{13}$ & 4&10&15&30& $W^{30}$, $W^6X^6$, $X^5Y$, $Y^3$, $Z^2$\cr
$Z_{11}$ & 6&8&15&30& $W^{30}$, $W^6Y^3$, $X^5$, $XY^3$, $Z^2$\cr
\noalign{\hrule}
$E_{14}$ & 3&8&12&24& $W^{24}$, $W^6X^6$, $X^4Z$, $Y^3$, $Z^2$\cr
$Q_{10}$ & 6&8&9&24& $W^{24}$, $W^6Z^2$, $X^4$, $Y^3$, $XZ^2$\cr
\noalign{\hrule}
$Z_{12}$ & 4&6&11&22& $W^{22}$, $W^6X^4$, $W^4Y^3$, $X^4Y$, $XY^3$, $Z^2$\cr
\noalign{\hrule}
$W_{12}$ & 4&5&10&20& $W^{20}$, $W^{10}Y^2$, $W^2X^2Y^2$, $X^5$, $Y^2Z$,
$Z^2$\cr
\noalign{\hrule}
$Z_{13}$ & 3&5&9&18& $W^{18}$, $W^6X^4$, $W^3Y^3$, $X^3Z$, $XY^3$, $Z^2$\cr
$Q_{11}$ & 4&6&7&18& $W^{18}$, $W^6X^3$, $W^4Z^2$, $X^3Y$, $Y^3$, $XZ^2$\cr
\noalign{\hrule}
$W_{13}$ & 3&4&8&16& $W^{16}$, $W^4X^4$, $W^4Y^3$, $X^4Y$, $Y^4$, $Z^2$\cr
$S_{11}$ & 4&5&6&16& $W^{16}$, $W^6Y^2$, $W^4Z^2$, $X^4$, $XZ^2$, $Y^2Z$\cr
\noalign{\hrule}
$Q_{12}$ & 3&5&6&15& $W^{15}$, $W^6X^3$, $W^3Z^2$, $X^3Z$, $XZ^2$, $Y^3$\cr
\noalign{\hrule}
$S_{12}$ & 3&4&5&13& $W^{13}$, $W^4X^3$, $W^3Z^2$, $WY^3$, $X^3Y$, $XZ^2$,
$Y^2Z$\cr
\noalign{\hrule}
$U_{12}$ & 3&4&4&12& $W^{12}$, $W^4Y^2$, $W^4Z^2$, $X^4$, $Y^2Z$, $YZ^2$\cr
\noalign{\hrule}
}}\qed$$

\enddemo

\remark{Remark 4.3.10}
One should choose the polyhedron rather carefully,
since the restricted K\"ahler class $\omega$ of
a general K\"ahler class of $\tilde \P_\Delta$
is sometimes not general in $(S_X)_\R$.
Namely, there might be some non-zero (non-effective) integral
algebraic 2-cocycle $\alpha$ such that $\omega\cup\alpha = 0$.
This is where $L_0$ appears.
Thus it is necessary to resort to a subpolytope of the full Newton polytope.
\endremark

\chap{4.4. Examples.}
We give some illuminating examples other than Arnold's singularities.

\example{Example 4.4.1}
Let us think about a pair of dual weight systems
$W_a = (2,3,6;12)$ and $W_b = (2,4,5;12)$.
In these cases $a_0 = b_0 = 1$.
We will denote the homogeneous coordinates by
$W=X_0$, $X=X_1$, $Y=X_2$ and $Z=X_3$ and
the inhomogeneous ones by small letters.
A polytope is represented by its vertex set $\left< \text{vertices} \right>$.

The weighted magic square $C$ is
$\left(\matrix 0 & 2 & 1\\ 3 & 2 & 0\\ 0 & 0 & 2\endmatrix\right)$
whose determinant is $(-12)$.
Thus
$B =
\left(\matrix -1 & 1 & 0\\2 & 1 & -1\\-1 & -1 & 1\endmatrix\right)$,
$-B^{-1} =
\left(\matrix 0 & 1 & 1\\1 & 1 & 1\\1 & 2 & 3\endmatrix\right)$.

We use $(\alpha_1,\dots,\alpha_n)$
for the coordinates coming from $M(W_a) = \Z^{\oplus n}$
and $[\beta_1,\dots,\beta_n]$ for those from $M(W_b)$.
Then $[\beta_1,\dots,\beta_n]B=(\alpha_1,\dots,\alpha_n)$ holds
for the vertices in $M_\Q$ and the hyperplanes in $N_\Q$.

$\overline\Delta = \left<W^{12}, X^6, Y^4, Z^2\right>$.
The coordinates of the vertices are :
$(-1,-1,-1)$, $(5,-1,-1)$, $(-1,3,-1)$, $(-1,-1,1)$;
$[-2,-3,-5]$, $[-2,2,1]$, $[2,0,-1]$, $[0,0,1]$,
respectively.
This is a reflexive simplex.
In fact, the dual is given by
$(\overline\Delta)^* = \left<W^{12}, Y^3,X^2Y^2,XZ^2\right>$.

Let $X_0$ be the singularity $\{x^6+y^4+z^2=0\}$ in $\C^3$.
This is an isolated hypersurface singularity with Milnor number
$\mu = 15$.
In the table of Arnold\cite{AGV},
$X_0$ belongs to class $W_{1,0}$.

$Z = \{x^6 + y^4 + z^2 + 1 = 0\}$ is a smooth surface in $\C^3$ thus
a Milnor fiber of $X_0$.

Ebeling \cite{E} has calculated Milnor lattices of
isolated bimodal singularities.
The Milnor lattice $H_2(Z)$ is a direct sum of $U$ and \par
\noindent\hphantom{o-}o\hphantom{-o-o-o-o-}o\par\vskip-8pt
\noindent\hphantom{p-}$\shortmid$\hphantom{-O-o-o-o-}$\shortmid$\par\vskip-8pt
\noindent o-o-o-o-o-o-o-o-o-o-o,
where o is a vanishing cycle.

The resolution graph of $X \to \bar Z$ is \par
\noindent\hphantom{o-o-}o $c$\par\vskip-8pt
\noindent\hphantom{p-o-}$\shortmid$\par\vskip-8pt
\noindent o-o-o-o-o\par\vskip-8pt
\noindent \hphantom{p-o-}$\shortmid$\par\vskip-8pt
\noindent $a$ $b$ o $b'\,a'$\par\vskip-4pt
\noindent \hphantom{o-o-}$c'$\par
\noindent
where o represents a smooth rational curve with self-intersection $(-2)$
and the central curve is $\{W=0\}$.
$L_D$ is generated by the central curve and
$(-4)$-elements $a+a'$, $b+b'$ and $c+c'$.
$L_0$ is generated by $a-a'$, $b-b'$ and $c-c'$.

Next think about the dual polygon $\overline\Delta^*$.
The singularity $\{y^3 + x^2y^2 + xz^2 = 0\}$ has a singular curve
$\{x = 0\}$.
We have a Milnor fiber $\{y^3 + x^2y^2 + xz^2 + 1 = 0\}$ and
its weighted compactification $\bar Z = \{Y^3 + X^2Y^2 + XZ^2 + W^{12} = 0 \}$
in $\P(1,2,4,5)$.
In this case, the restriction of the divisor $\{X = 0\}$ decomposes as
a sum of three divisors $E_i : \{X = Y+\zeta^iW^4 = 0\}$ ($i=0,1,2$)
where $\zeta$ is a nontrivial cubic root of unity.
The other intersections with singular 2-dimensional strata of $\P$
are irreducible.

The singular set of $\bar Z$ consists of one $D_8$ at $\{W=Y=Z=0\}$,
one $A_4$ at $\{W=X=Y=0\}$ and one $A_1$ at $\{W=X^2+Y=Z=0\}$.
The resolution graph at infinity (including the divisor $W=0$) becomes\par
\noindent\hphantom{}$a'$o\hphantom{-o-o-o-o-o-}o\par\vskip-8pt
\noindent\hphantom{p-}$\shortmid$\hphantom{-O-o-o-o-o-}$\shortmid$\par\vskip-8pt
\noindent o-o-o-o-o-o-o-o-o-o-o-o.\par\vskip-4pt
\noindent $a$ $b$ $c$\par
This graph is one longer than that of the $L_G(\overline\Delta)$.
$L_0$ is generated by $a-a'$, $E_0-E_1$ and $E_1-E_2$.
$L_D$ is generated by $a+a'$ and other o's.
Change the base as $a+a'\to a'+a+b$ (other o's are same),
we have the same graph as the $L_G(\overline\Delta)$.

Let us take other polygons.
In our duality, every reflexive polygon in $M_\Q$ should contain
${\overline\nabla}^* = <W^{12}, X^3Y^2, Y^2Z, Z^2>$,
or rather $\nabla^* = <W^{12}, Y^2W^6, X^3Y^2, Y^2Z, Z^2>$,
where $\nabla$ is the maximum Newton polygon in $N$,
that is the convex hull of $\overline\nabla \cap N$.

Take
$\Delta = <W^{12}, W^3Y^3, W^2X^5, WXY^3, X^3Y^2, X^3Z, Y^2Z, Z^2>$,
which has the same $L_G$ as $\overline\Delta$ and $L_0(\Delta) = 0$.
Thus $\Delta$ represents a Milnor fiber of the singularity of class $W_{1,0}$
and the dual graph of $L_D$ coincides with the resolution graph as above.
Its dual $\Delta^*$ is
$<W^{12}, W^7Z, W^6XY, W^5XZ, W^3YZ, X^2Y^2, Y^3, XZ^2>$.
\endexample

\example{Example 4.4.2}
The weight system $W = (2,3,4;10)$ is self-dual.
Take the following six reflexive polygons:
$$\align
\Delta_1 &= <W^{10}, W^2X^4, X^3Z, X^2Y^2, XZ^2, Y^2Z, W^2Z^2, WY^3>, \\
\Delta_2 &= <W^{10}, W^4X^3, WX^3Y, X^2Y^2, XZ^2, Y^2Z, W^2Z^2, WY^3>,\\
\Delta_3 &= <W^{10}, W^6X^2, W^2X^2Z, X^2Y^2, XZ^2, Y^2Z, W^2Z^2, WY^3>,\\
\Delta_4 &= <W^{10}, W^6X^2, W^2X^2Z, X^2Y^2, XZ^2, Y^2Z, W^6Z, W^4Y^2>,\\
\Delta_5 &= <W^{10}, W^8X, W^3X^2Y, X^2Y^2, XZ^2, Y^2Z, W^6Z, W^4Y^2>
\text{and}\\
\Delta_6 &= <W^{10}, W^4XZ, X^2Y^2, XZ^2, Y^2Z, W^6Z, W^4Y^2>.
\endalign$$
These polygons all satisfy $L_0 = 0$ and
$\Delta_i^* \cong \Delta_{7-i}$.
$\Delta_1$ corresponds to the Milnor fiber
of the singularity of class $S_{1,0}$.

The dual graphs of $L_D(\Delta_i)$ is as follows:\par
{%\eightpoint
$\Delta_1$\par
\noindent\hphantom{o-}o\par\vskip-8pt
\noindent\hphantom{p-}$\shortmid$\par\vskip-8pt
\noindent\hphantom{o-}o\par\vskip-8pt
\noindent\hphantom{p-}$\shortmid$\par\vskip-8pt
\noindent o-o-o-o-o.\par\vskip-8pt
\noindent \hphantom{p-}$\shortmid$\par\vskip-8pt
\noindent \hphantom{o-}o\par
$\Delta_2$\par
\noindent\hphantom{o-o-}o\par\vskip-8pt
\noindent\hphantom{p-o-}$\shortmid$\par\vskip-8pt
\noindent\hphantom{o-}o\hphantom{-}o\par\vskip-8pt
\noindent\hphantom{p-}$\shortmid$\hphantom{p}$\shortmid$\par\vskip-8pt
\noindent o-o-o-o-o-o.\par
$\Delta_3$ and $\Delta_4$\par
\noindent\hphantom{o-o-o-}o\par\vskip-8pt
\noindent\hphantom{p-o-o-}$\shortmid$\par\vskip-8pt
\noindent\hphantom{o-}o\hphantom{-o-}o\par\vskip-8pt
\noindent\hphantom{p-}$\shortmid$\hphantom{-O-}$\shortmid$\par\vskip-8pt
\noindent o-o-o-o-o-o-o.\par
$\Delta_5$\par
\noindent\hphantom{o-o-o-o-}o\par\vskip-8pt
\noindent\hphantom{p-o-o-o-}$\shortmid$\par\vskip-8pt
\noindent\hphantom{o-}o\hphantom{-o-o-}o\par\vskip-8pt
\noindent\hphantom{p-}$\shortmid$\hphantom{-O-o-}$\shortmid$\par\vskip-8pt
\noindent o-o-o-o-o-o-o-o.\par
$\Delta_6$\par
\noindent\hphantom{o-o-o-o-o-}o\par\vskip-8pt
\noindent\hphantom{p-o-o-o-o-}$\shortmid$\par\vskip-8pt
\noindent\hphantom{o-}o\hphantom{-o-o-o-}o\par\vskip-8pt
\noindent\hphantom{p-}$\shortmid$\hphantom{-O-o-o-}$\shortmid$\par\vskip-8pt
\noindent o-o-o-o-o-o-o-o-o.\par
}

There are inclusions of reflexive polytopes
$\Delta_i \supset \Delta_j$ for $i < j$
and thus specializations of the families of weighted hypersurfaces.
As the family specializes,
$2A_1$ singularities at infinity get `closer'
to $A_3$, $D_4$, $D_5$ and finally to $D_6$ singularities.
\endexample

\enddocument